**Title:**

Behavioural predictors of math anxiety

**Authors:**

Name: Chen Mun Yip Kinnard
Institution: Nanyang Technological University
Email: chen1308@e.ntu.edu.sg

Name: Azilawati Jamaludin
Institution: National Institute of Education, Nanyang Technological University
Email: azilawati.j@nie.edu.sg

Name: Aik Lim Tan
Institution: National Institute of Education, Singapore
Email: aiklim.tan@nie.edu.sg
ORCiD: https://orcid.org/0000-0001-5910-7003

**Corresponding author:**

Name: Aik Lim, TAN
Institution: National Institute of Education, Singapore (1 Nanyang Walk, Singapore 637616)
Email: aiklim.tan@nie.edu.sg

**Declarations of Interest:** The authors declare no conflict of interest.

**Funding**

This research was funded by the Singapore National Research Foundation (NRF) under the Science of Learning Initiative (NRF2016-SOL002-003).

# Behavioural Predictors of Math Anxiety

*Abstract*

Math anxiety (MA) is a highly prevalent problem in education that has consistently shown to lead to poorer math performance. This study sought to investigate whether certain behaviours are predictive of MA among students. The study involved elementary school students (*n=124*) who were low-progressing in math and is part of an educational intervention program to improve their mathematical abilities through a neural-informed math game. Ten classification types of behavioural indicators were identified, such as counting out loud. A multiple linear regression was conducted where students' anxiety scores were regressed on these behavioural observations, along with their gender, mood, and math profile. Three behavioural observations were positively and significantly associated with their math anxiety. Implications and limitations of the study are discussed.

## INTRODUCTION

Math anxiety (MA) refers to a state of fear and apprehension when one is engaging with math (Ashcraft & Krause, 2007). It is a highly prevalent problem in many educational settings and persists across many different levels of education (Gunderson et al., 2018; Jain & Dowson, 2009), beginning in primary school (Harari et al., 2013) and even up until university education (Røykenes & Larsen, 2010). Furthermore, MA has been shown to lead to poorer math performance. For instance, a meta-analysis by Zhang et al. (2019) found that MA was robustly and negatively associated with math performance. In other words, higher MA tend to lead students to perform more poorly in math.

Researchers have also found that MA can have an adverse impact on one's initial learning of mathematics and subsequently, a poorer acquisition of math skills and advanced mathematical concepts as well (Krinzinger et al., 2009). For instance, MA has been shown to adversely affect basic numerical operations such as simple addition and subtraction – operations that serve as the foundation for more complex mathematical abstraction (Ashcraft, 2002; Maloney et al., 2010). Indeed, Ashcraft and Moore (2009) found that individuals with MA tend to find it more difficult to engage in complex mathematical calculations.

It should also be noted that MA can both be the cause and the result of poorer math performance among students (Foley et al., 2017). This means that MA does not only lead to poorer math performance, but it can also occur due to having performed poorly in a standardized math test (Ashcraft & Krause, 2007). In the long-term, this can result in a feedback loop whereby after having performed poorly in a math test, students get MA and start to avoid math lessons (Ashcraft & Moore, 2009). Subsequently, they might also develop negative beliefs about their math abilities, which further exacerbates their MA and avoidance of math classes (Lent et al., 1991). For example, MA is often regarded as one of the primary reasons for the low enrolment in math courses among female students (Meece et al., 1982). This avoidance can also extend to mathematics-related professions, which further limits one's career opportunities (Eispino et al., 2017).

Similarly, a recent meta-analysis by Barroso et al. (2021) found that MA was negatively associated with math achievement, and that this negative association was present from childhood through to adulthood. It should be noted that poorer math performance also lead to negative outcomes in adulthood such as poorer financial planning (McKenna and Nickols, 1988) and lower confidence among math educators (Olango and Assefa, 2013).

Given that MA can lead to poorer math learning and performance both in the short-term and long-term, it is important to help students manage or ameliorate their MA. A systematic review by Balt et al. (2022) suggests that interventions that help manage one's MA also serve to improve one's math performance. It is therefore important for the early identification of students with higher levels of MA, in order to help them reduce their feelings of anxiety, and subsequently, improve their math performance.

Presently, the predominant way of identifying MA in students is through self-report measures (Carey et al., 2017; Luttenberger et al., 2018). However, this method of identifying MA among students in actual educational settings might not be very practical as time and manpower are needed to administer and collate the results of these measures (Wu et al., 2012). Furthermore, given that one's level of MA can change over time (Ashcraft & Krause, 2007), there might be a need for regular testing to obtain a more accurate measure of the students' MA. This suggests that there is a need for a more cost-effective way of identifying MA among students and across time.

One way of identifying MA might be through the identification of behaviours that are symptomatic of one's underlying MA. Anxiety can manifest itself in a variety of behaviours such as fidgeting and restlessness (Yadav et al., 2017). This suggests that it might be possible to identify MA in students through their behaviour when engaging with math-related tasks. This is important as any behaviour that is predictive of MA can help educators to easily identify students with high MA and, consequently help reduce their MA by either managing their emotions or improving their math performance (Balt et al., 2022). Thus, this study sought to investigate whether certain behaviours are predictive of MA among students.

## LITERATURE REVIEW

### Math Anxiety

The concept of MA was first introduced as a topic of study in academic research by Dreger and Alken (1957) to describe students' anxiety toward numbers. Presently, MA commonly refers to a state of fear and apprehension when one is engaging with math (Ashcraft & Krause, 2007) and is regarded as a primarily emotional response (Mammarella et al., 2015). It should be acknowledged that although most studies conceptualize MA as a single factor, some researchers suggest that MA consists of two dimensions: cognitive and affective (Wigfield & Meece, 1988). The cognitive dimension broadly refers to one's thoughts and concerns about math performance while the affective dimension includes one's emotions like nervousness regarding math testing (Wigfield & Meece, 1988).

MA can also be differentiated from test anxiety or general anxiety (Ashcraft & Ridley, 2005) For instance, different measures of MA are more highly correlated with one another as compared to test and general anxiety (Ashcraft & Ridley, 2005). In other words, MA has been shown to be sufficiently different from other forms of test or general anxiety, and can thus be considered to be an independent construct (Dowker et al., 2016).

There are many factors that can affect one's level of MA and they can be differentiated into three main categories: environmental, cognitive and personal factors (Rubinsten & Tannock, 2010). Environmental factors include the experiences of math learning in the classroom or at home. For instance, Mutodi and Ngirande (2014) found that negative experiences of math learning in both at home and in the classroom can lead to MA. Similarly, stressful teaching strategies like time testing (Ashcraft & Moore, 2009) and assigning math-related tasks as a form of punishment have been found to increase one's MA (Oberlin, 1982). Conversely, McNaught and Grouws (2007) suggests that teachers and parents should focus on creating a more positive environment for students when they are learning mathematics. They also suggested that parents should focus on exploring any fear and anxieties that the child might have towards learning mathematics, especially during the earlier stages of math learning. Furthermore, parents' own MA can have a negative intergenerational effect on their children's experiences of math learning, leading them to have poorer math learning and higher levels of MA (Maloney et al., 2015). The culture in which one is brought up could also influence one's level of MA, whereby students in countries with a greater emphasis on doing well in examinations could experience higher levels of MA (Dowker et al., 2016; Tan & Yates, 2011).

Personal factors include low self-efficacy (Kesici & Erdoğan, 2010; Maloney et al., 2011), low self-esteem, previous negative experiences (Mutodi & Ngirande, 2014; Rubinsten & Tannock, 2010) and even genetic vulnerabilities (Wang et al., 2014). For instance, Wang et al. (2014) found that one's genetic vulnerability to general anxiety also increases one's likelihood of developing MA. Cognitive factors involve having low intelligence and poor cognitive abilities in mathematics (Rubinsten & Tannock, 2010).

### Assessments of Math Anxiety

Math anxiety is assessed almost exclusively with self-report questionnaires, across all age groups, in both educational and research settings (Carey et al., 2017; Luttenberger et al., 2018). For adults, common measures include the Mathematics Anxiety Rating Scale (MARS; Richardson & Suinn, 1972) and the Revised Mathematics Anxiety Rating Scale (R-MARS; Baloğlu & Zelhart, 2007). There are also questionnaires developed for younger populations such as primary school children, involving more pictorial rating scales and greater focus on more concrete math situations. Examples include the Scale for Early Mathematics Anxiety (SEMA; Wu et al., 2012), and the Children's Attitude to Math Scale (James, 2013).

As with most self-report measures, the accuracy of these questionnaires are subjected to the participants' truthfulness and the accuracy of their own self-perceptions (Dowker et al., 2016). To counteract this problem, some studies utilized physiological measures to determine the individual's level of MA, such as the level of cortisol secretion (Mattarella-Micke et al., 2011). Neurological measures such as EEG recordings (Núñez-Peña & Suárez-Pellicioni, 2015) and functional MRI scans have also been used to measure and map out MA (Pletzer et al., 2015). However, such measures might not be very practical in actual educational settings as time and manpower are needed to administer and collate the results of these measures (Pletzer et al., 2015; Wu et al., 2012). Given that one's level of MA can change over time (Ashcraft & Krause, 2007), regular testing might be required to obtain a more accurate measure of students' MA. This suggests that there is a need for a more cost-effective and adaptive way of measuring one's level of MA among students and across time.

**Behavioural Predictor of Math Anxiety**

One possible indicator of MA could be the behaviour of the individual while he or she is engaging with math-related tasks. MA can affect the individual and manifest itself at various levels, such as the physiological level (Hunt et al., 2017), the psychological or cognitive level (Hunt et al., 2014), and the emotional level (Papousek et al., 2012).

MA can also affect individuals on a behavioural level (Pizzie & Kraemer, 2017). For instance, individuals who have higher levels of MA tend to display higher instances of disengagement, exemplified by the avoidance of mathematical stimuli such as math classes or math-related work (Preis & Biggs, 2001; Pizzie & Kraemer, 2017). Similarly, some researchers assert that MA invokes three types of reactions: affective, cognitive, and behavioural responses (Olango, 2016; Ashcraft, 2019). Olango (2016) states that these behavioural responses fall on a continuum with approach and avoidance behaviours as the continuum extremes. Olango (2016) also describes how the three domains can interact with one another in a reinforcing way. This could produce either a negative cycle, whereby avoidance behaviour reciprocally interacts feelings of anxiety and worry, or whereby one's approach behaviour reciprocally interacts with thoughts of math achievement and feelings of confidence, creating a positive and self-reinforcing cycle.

Similarly, Pries and Biggs (2001) outlines how the negative cycle of avoidance of mathematical stimuli occurs in four phases and how it is self-reinforcing. In the first phase, the individual has a negative experience with mathematics or mathematics situations. This leads to the second phase, whereby the individual starts to avoid mathematics due to the negative experience previously. For the third phase, after avoiding mathematics, the individual has less opportunities to improve their mathematics skills, consequently becoming less proficient in and prepared for mathematics (Dowker et al., 2016; Pries & Biggs, 2001). In the fourth phase, as a consequence of being less prepared for mathematics, the individual has poorer mathematics performance. This creates a negative experience with mathematics, leading the individual back to phase one and developing a vicious cycle of mathematics avoidance and MA (Pries & Biggs, 2001; Ramirez et al., 2018). Such avoidance behavior is so common that Ashcraft and Moore (2009) describe it as "an overriding characteristic of math-anxious individuals" (p. 201).

This suggests that certain behaviours could be a manifestation of one's MA and consequently be used to predict the presence of it. This is consistent with the symptomology of anxiety in clinical research, whereby certain behaviours have been shown to be symptomatic of one's anxiety (American Psychiatric Association, 2013; Beck et al., 1988; Spitzer et al., 2006). For instance, the Beck Anxiety Inventory (BAI) – a self-report measure of general anxiety symptomatology – indicates that certain behaviours such as being "shaky", "hands trembling" and "wobbliness in legs" are common symptoms of anxiety (Beck et al., 1988). Similarly, the Generalised Anxiety Disorder-7 scale

(GAD-7) by Spitzer et al. (2006) – a common self-report tool to determine the severity of anxiety symptoms – includes the behavioural symptom of "being so restless that it is hard to sit still". The clinical literature therefore suggests that the presence of certain behaviours is predictive of one's anxiety, signifying that certain behaviours might also be predictive of one's MA.

**METHOD**

**Participants**

A total of 124 Singapore primary school students – aged 6 to 7 years old – were recruited from primary schools in Singapore for this study. Of the 124 students, 48.39% were females (n = 60) and 51.61% were males (n = 64). The students had different math profiles; they were labelled as either Typically Developing (TD; 15.32%; n = 19), Low-Progress Learners (LP; 41.13%; n = 51), or at-risk of Developmental Dyscalculia (DD; 43.55%; n = 54). The students were categorised into the three groups based on a national numeracy screening assessment administered by the schools when students first enter primary school. The students were grouped according to four levels based on their scores, ranging from 0 to 3. Students with scores 2 or below are considered low-progressing and will need to attend a math pull-out program aimed at providing more support for these students. Students who score a '3' are considered Typically Developing. In this particular study, the students also completed the Dyscalculia Screener (Butterworth, 2003), which identified students who may be at-risk of developmental dyscalculia.

The study was conducted in accordance with the Declaration of Helsinki, and the protocol was approved by the NTU Institutional Review Board (Reference Number: IRB-2017-10-030).

**Measures**

*Math Anxiety*

To measure the students' math anxiety, the Scale for Early Mathematics Anxiety (SEMA) by Wu et al. (2012) was used. The SEMA scale is a 20-item self-report measure that requires participants to indicate how anxious they would feel in certain situations involving mathematics. It utilizes a 5-point Likert scale, ranging from 1 (not nervous at all) to 5 (very, very nervous). Each of the 5 points on this Likert scale are accompanied by a drawing of an anxious face – each with different gradations of anxiety based on the number (or the intensity of the anxiety) they represent. The accompaniment of the anxious faces helps the children to better identify and choose the number that best represent their levels of anxiety (Wu et al., 2012).

The measure can be differentiated into its two factors: 1) Numerical Processing Anxiety (NPA) factor and 2) Situational and Performance Anxiety (SPA) factor. The NPA items require respondents to indicate how anxious they would feel if they had to solve certain math questions (e.g., "Is this right? 9 + 7 = 18"), while the SPA factor focuses on situations where mathematics is involved (e.g., "You are in math class and your teacher is about to teach something new"). Each factor consists of 10 items.

The total SEMA score is calculated by summing all 20 item scores while the total score of its factors (e.g., NPA and SPA) is calculated by summing the scores of their corresponding 10 items. Higher scores represent higher levels of MA.

*Math Game – 'Number Beads'*

Number Beads is a digital game that was developed by Butterworth & Laurillard (2017), that sought to incorporate and apply the findings from neuroscience and cognitive science research to help low attaining learners (e.g., individuals with dyscalculia) in learning mathematics.

In the game, there is a target set of beads at the top of the screen (grouped with the word 'Target'; refer to Figure 1), with the main objective of the game being to construct the target set by either splitting or joining the different colour-coded beads available in the play area. More detailed information about the game 'Number Beads' can be found in Laurillard (2016).

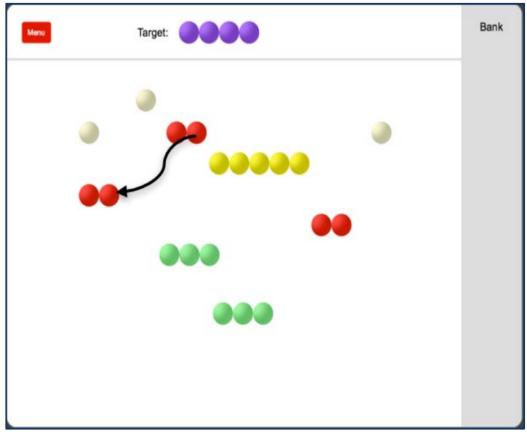

*Figure 1.* A screenshot of the Number Beads math game from Laurillard (2016).

*Behavioural Observations*

Given limited research on behaviours predictive of MA, this study took a thematic analysis approach towards the formulation of the different categories of behaviours that would be investigated. More specifically, the researchers observing the students were instructed to write down any notable behaviours that the students exhibited while engaging in the math game.

Based on these notable behaviours, they were then clustered and categorised into 10 different behavioural observations (BO; e.g., BO1 to BO10):

1. **Counting out loud (BO1)**

This involves counting by saying the numbers out loud one by one (e.g., "three, four, five…").

2. **Verbalization of thought process (BO2)**

This refers to the students' verbalization of their thought process as to how they would carry out the arithmetical calculations (e.g., "I can get 6 by taking 2 from 8").

3. **Utterance of other comments or sounds (BO3)**

Any other verbal utterances that do not fall into the above two behavioural observations would be included in this behavioural observation. Examples include humming, singing of songs, making nonsensical sounds or utterances that bear no meaning or relation to the arithmetical tasks at hand.

4. **Use of fingers to assist in counting (BO4)**

This includes using fingers to point to the 'beads' on the computer screen while counting and/or using one's fingers as representations of digits to assist in arithmetic calculations.

5. **Checking on and interacting with peers (BO5)**

This behavioural observation broadly refers to any of the students' interactions with their peers. Examples include looking at their classmates, asking how they are doing and/or competing with one another to achieve the highest scores.

6. **Gross hand and leg movements (BO6)**

This involves any gross hand and leg movements such as fidgeting, shaking of one's legs or playing with stationery.

7. **Asking questions when in doubt (BO7)**

Any act of clarification regarding how to play the math game or asking for assistance to complete a particular set by the students would be included in this behavioural category.

   8. **Celebrating success (BO8)**

This involves any behaviour that signifies a celebration of having successfully completed a level or a task. Examples of this behavioural observation include punching the air and saying "Yes! I did it!" or raising both hands in the air and cheering.

   9. **Looking elsewhere or being distracted (BO9)**

This refers to any behaviour that is suggestive of a lack of attention or being distracted such as looking elsewhere. If the student is looking at their peers, it will not fall under this behavioural observation and will instead fall under the behavioural observation of 'Checking on and interacting with peers'.

   10. **Random splitting and joining of Number Beads (BO10)**

This refers to the beads in the math game that participants manipulate (e.g., separating a group of 3 beads from a larger group of 8 beads to do the calculation for '8-3=5'). Even though this is not strictly a physical behaviour, there is a strong parallel between this behaviour and the actual manipulation of physical objects or what is termed as 'manipulatives' to do arithmetic (Jones & Tiller, 2017).

*Mood & Math Profile*

Based on the observations by the research team, each student's mood while engaging in the math game was labelled as one of four moods: 1) happy, 2) neutral, 3) bored or 4) stressed.

The math profile of the students – as either TD, LP or DD – was determined by both school's criteria (using the Ministry of Education's Early Numeracy Indicator) and the Dyscalculia Screener (Butterworth, 2003).

**Procedure**

The students were first instructed to complete the SEMA measure. Aligned with Wu et al. (2012), a trained researcher assisted the students in the completion of this measure by helping to read the questions aloud and asking the students to point to the face that best represented how anxious they felt. The students were also encouraged to ask questions and clarify the meaning of the questions if they did not understand the scale's items.

Students were subsequently invited to play the neural-informed math game – Number Beads – on a designated laptop. During the game, the researchers took note of the students' mood and any notable behaviours that they exhibited.

The administration of SEMA and game play took place in school on separate sessions where the game was played for about 2-4 hours per student.

**Data Analysis**

Firstly, descriptive analyses of the following variables were conducted: 1) mood; and 2) SEMA scores.

Secondly, the Chi-square test was employed to investigate whether gender, affect and the behavioural observations were significantly related to the students' math profiles.

Thirdly, multiple linear regression analyses were conducted where the SEMA scores were regressed on the 1) behavioural observations, and the students' 2) gender, 3) mood and 4) math profile. Given that there is no research on such behavioural indicators of math anxiety among Singapore learners, a stepwise regression was chosen as it would allow for the selection of a model that provides the most efficient prediction of MA (Aiken & West, 1991). For the stepwise regression analysis, the PIN will be set at .05 and POUT will be set at .10.

As SEMA scale can be broken down into two factors: NPA & SPA factors. Thus, the stepwise regression analyses were conducted thrice; each with a different dependent variable (DV). For the first model, the DV is the total SEMA score. The second model's DV is the total NPA score, while the third model's DV is the total SPA score.

**RESULTS**

**Descriptive Analyses**

*Mood*

Based on the researchers' observations, 51.61% of the students were happy (n = 64), 39.52% were of neutral mood (n = 49), 8.06% were bored (n = 10) and 0.81% was stressed (n = 1).

*SEMA scores*

With regards to the total SEMA score, the mean score across all participants was 23.5 (possible range of 0 – 80), with a standard deviation (SD) of 15.91. For each item, the average score was 1.18 (out of 4), with an SD of 0.80.

The SEMA scores can be differentiated into their two factors: the NPA factor and the SPA factor. For the NPA factor, the mean score across all participants was 13.24 (possible range of 0 – 40), with a SD of 8.60. For each item within the NPA factor, the average score was 1.32 (out of 4), with an SD of 0.86. For the SPA factor, the mean score across all participants was 10.26 (possible range of 0 – 40), with a SD of 8.76. For each item within the SPA factor, the average score was 1.03 (out of 4), with an SD of 0.88.

**Math Profile and Gender, Affect & Behavioral Observations**

The chi-square tests revealed that affect, ($X^2$ (2, *N* = 124) = 5.696, *p* = .458), and gender were not significantly related to math profile ($X^2$ (2, *N* = 124) = 0.611, *p* = .737).

Math profile was only significantly related to one of the behavioral observations, namely BO5 - Checking on and interacting with peers ($X^2$ (2, *N* = 124) = 8.540, *p* = .014). More specifically, those with development dyscalculia were more likely to exhibit the behavior of 'checking on and interacting with peers' when engaging with math, as compared to those who were 'low-progress learners' and 'typically developing'. It was a small to medium effect size (*V* = .262; Cohen, 1988).

**First regression model: Total SEMA score**

The first regression model regressed the students' 1) affect, 2) gender, 3) math profile and 4) behavioural observations on the total SEMA scores, using a stepwise approach.

The first predictor added into the regression model was BO1 – Counting out loud. It had the largest coefficient of determination ($R^2$ = .061, F(1, 122) = 7.863, *p* < 0.01). After BO1 was present in the model, BO10 – Random splitting and joining of Number Beads – had the next biggest change in the coefficient of determination ($R^2$ = .050, F(1, 121) = 6.825, *p* < 0.05), and was thus added into the model.

The concluding model revealed BO1 and BO10 to be significant factors in affecting total SEMA scores (refer to Table 1). This set of predictors explained 11.1% of the variance in total SEMA scores ($R^2$ = .111, F(2, 121) = 7.531, *p* < 0.001).

Table 1. Summary of the Multiple Regression Analysis for the First Model on Total SEMA scores.

| Variable | Unstandardised Coefficients | | Standardised Coefficients | t | Sig. |
|---|---|---|---|---|---|
| | B | Std. Error | β | | |
| Constant | 19.078 | 1.777 | | 10.737 | .000 |
| BO1 – Counting out loud | 7.807 | 2.934 | .229 | 2.660 | .009 |
| B010 – Random splitting and joining of Number Beads | 8.409 | 3.219 | .225 | 2.612 | .010 |

**Second regression model: Total NPA score**

Similarly, the second regression model regressed the students' 1) affect, 2) gender, 3) math profile and 4) behavioural observations on the total NPA scores, using a stepwise approach.

The behavioural observation – BO1 (Counting out loud) – was the only predictor added in this regression model (R2 = .047, F(1, 122) = 5.983, p < 0.05). In other words, this model showed that BO1 was the only factor that was significantly associated with the total NPA scores (refer to Table 2). This model explained 4.7% of the variance in total SEMA scores (R2 = .047, F(1, 122) = 5.983, p < 0.05).

Table 2. Summary of the Multiple Regression Analysis for the Second Model on Total NPA scores.

| Variable | Unstandardised Coefficients | | Standardised Coefficients | t | Sig. |
|---|---|---|---|---|---|
| | B | Std. Error | β | | |
| Constant | 11.988 | .914 | | 13.117 | .000 |
| BO1 – Counting out loud | 3.986 | 1.630 | .216 | 2.446 | .016 |

**Third regression model: Total SPA score**

Lastly, the third regression model regressed the same four classes of variables on the total SPA scores, using a stepwise approach.

The first predictor added into the regression model was BO10 – Random splitting and joining of Number Beads ($R^2$ = .075, $F$(1, 122) = 9.906, $p$ < 0.01). BO1 – Counting out loud – was then subsequently added into the model ($R^2$ = .046, $F$(1, 121) = 6.319, $p$ < 0.05). Lastly, the model then included BO9 – Looking elsewhere or being distracted ($R^2$ = .033, $F$(1, 120) = 4.739, $p$ < 0.05).

The concluding model showed that BO1, BO9 and BO10 were significant factors in affecting total SPA scores (refer to Table 3). This set of predictors explained 15.4% of the variance in total SEMA scores ($R^2$ = .154, $F$(1, 120) = 7.304, $p$ < 0.001).

Table 3. Summary of the Multiple Regression Analysis for the Third Model on Total SPA scores.

| Variable | Unstandardised Coefficients | | Standardised Coefficients | t | Sig. |
|---|---|---|---|---|---|
| | B | Std. Error | β | | |
| Constant | 6.818 | 1.049 | | 6.501 | .000 |
| BO10 – Random splitting and joining of Number Beads | 4.564 | 1.769 | .221 | 2.580 | .011 |
| B01 – Counting out loud | 4.419 | 1.592 | .235 | 2.776 | .006 |
| BO9 – Looking elsewhere or being distracted | 3.694 | 1.697 | .187 | 2.177 | .031 |

**Discussion**

This study sought to investigate whether affect, gender, math profile and certain behavioural observations are predictive of MA. This study concluded with 3 different regression models: 1) total SEMA scores, 2) total NPA scores and 3) total SPA scores.

Is gender related to MA?

The findings revealed that there was no significant relationship between gender and MA for all three regression models. This is consistent with many previous studies which found that there were no gender differences in MA for elementary school students (i.e., students aged 5 to 10; Erturan & Jansen, 2015; Harari et al., 2013; Kucian et al., 2018; Schleepen & Van Mier, 2016).

Is affect (mood) related to MA?

There was no significant association found between affect and MA for all three regression models. In other words, the findings suggest that one's affect is not related to one's level of MA. This insignificant association could be due to the students' differential abilities in emotional regulation. One form of emotional regulation is expressive suppression, which refers to the attempt to conceal one's feelings and physiological state (Cohen et al., 2021). Prior studies found that an individual's ability to engage in emotional regulation is linked to MA (Brooks, 2014; Pizzie & Kraemer, 2021). This suggests that the student's differential abilities in emotional regulation could be a moderating factor in the relationship between affect and MA, which might explain the lack of a significant association between these two variables.

Are there differential associations of MA with math learning profiles?

The results suggest that there is no significant relationship between math profile and MA for all regression models. That is, there was no significant differences in the MA scores among students with different math profiles. The lack of a significant relationship between math profile and MA in the current study could be due to differences in the testing contexts. Unlike previous studies which measured the students' MA in traditional educational contexts (e.g., before and after a standardized math test; Zhang et al., 2019), the current study measured the students' MA in their mathematics classrooms and prior to the gameplay.

Are there behavioural predictors of MA?

The findings revealed that three behavioural observations were significantly and positively associated with MA: namely, 1) the act of counting out loud, 2) randomly splitting and joining of the Number Beads, and 3) looking elsewhere or being distracted.

BO1 – counting out loud – was consistently included in all three regression models. This suggests that the behaviour of counting out loud could be a predictor of a student with MA. This finding seems to be consistent with the current literature. Past findings suggest that individuals with MA tend to have smaller work memory capacity when engaging in arithmetic tasks as they often have intrusive anxious thoughts that take up working memory resources (Ashcraft & Krause, 2007; Maloney et al., 2010; Shi & Liu, 2016). It was also found that one's working memory capacity mediates the relationship between MA and visual enumeration (Maloney et al., 2010). Given that saying material out loud has been shown to improve one's memory of it (Forrin & MacLeod, 2018), the act of counting out loud could similarly be a compensatory strategy to aid one's reduced working memory capacity.

BO10 – random splitting and joining of Number Beads was included in two regression models – where the predictors were regressed on SEMA scores and SPA scores. This suggests a positive relation between student's overall MA and his/her situational and performance anxiety. This is consistent with current literature on the negative relationship between MA and math performance (Zhang et al., 2019). This behavioural observation (BO) of the random management of the digital manipulatives suggests that the student may not know how to procedurally solve the arithmetic task, suggesting underlying math struggles manifested by poor math performance. Furthermore, MA can both be the cause and the result of poorer math performance (Foley et al., 2017). Taken together, BO10 could be symptomatic of one's poorer math performance and or MA. However, given that this behavioural observation was not significantly associated with NPA scores, it suggests that this behaviour only occurs for individuals who are anxious about being evaluated on their ability to do math problems, and not for those who feel anxious in nonevaluative situations.

BO9 – looking elsewhere or being distracted – was only included in one regression model, where it was focused on SPA scores. This suggests that this behaviour is positively related to one's situational and performance anxiety. This is consistent with prior studies on the behavioural responses that are commonly associated with math anxiety, such as behavioural disengagement (Pizzie & Kraemer, 2017). Similarly, Eysenck et al. (2007) also showed that anxiety decreases one's focus on the task at hand and makes one more vulnerable to distraction. An inability to focus on the task at hand and being distracted by other stimuli (e.g., one's intrusive thoughts) are common characteristics of those with high levels of anxiety (Ramirez et al., 2016). The behavioural observation of looking elsewhere or being distracted seems to be symptomatic of this tendency towards disengagement and distraction. However, since BO9 was not significantly related to SEMA and NPA scores, this suggests that this avoidance behaviour only occurs for individuals who feel anxious when they are being evaluated on their ability to do math problems, and not for those who feel anxious in nonevaluative situations.

**Theoretical and practical implications**

Given that there has been no study investigating the relationship between behavioural observations and MA particularly in the Singapore education context, this study helped to bridge this literature gap. The current study therefore provides researchers with preliminary insights into MA through an ethnographic lens, affording opportunities for future research to be done in this area. The findings provide some practical implications in the context of math education.

Instead of relying on self-report measures, educators can tap on their behavioural observations of classroom manifestations to anticipate quickly and cost-effectively socio-emotional issues related to math learning such as that of Math Anxiety. It would also be less onerous for educators to rely on such BO instead of having students complete interim self-report measures of MA. Such behavioural observations could also be used for educational contexts in which students are learning to count by manipulating physical objects (Jones & Tiller, 2017). The findings suggest that the random manipulation of physical objects (e.g., beads) could also be used as an indicator of one's MA. Such observations are therefore able to provide teachers with early indicators to identify students with possible MA, thereby enabling them to intervene early and provide the necessary support for these students.

**Limitations and Future research**

This study has two main limitations. Firstly, this study was conducted in the context of a math game. More specifically, the students that were being observed in the present study were playing a math game. Such educational games tend to be less anxiety-inducing than traditional educational settings due to the lower stakes and expectations (Marshall et al., 2016; Rocha & Dondio, 2021; Taylor & Mohr, 2001). Thus, the behaviours that were predictive of MA in the current study's context might not be generalizable to traditional educational settings. More research will be needed to ascertain the utility of such behavioural predictors in real-world educational settings.

Secondly, this study was conducted in Singapore, an East Asian country with a high level of math proficiency among its students (Programme for International Student Assessment, 2018). Many studies have shown that emotional suppression – which refers to the concealment of one's true emotions (Schouten et al., 2020) – is more common in Asian countries like Singapore and Japan as compared to Western countries like Belgium and the United States (Butler et al., 2007; Cheung & Park, 2010; Soto et al., 2011). This suggests that students from different countries may have differing levels of emotional suppression or regulatory abilities when engaging in math-related activities and consequently differing expressions of MA.

**CONCLUSION**

This study sought to investigate the relationship between behavioural observations and MA. The findings suggest that behavioural observations like counting out loud, randomly splitting and joining Number Beads (or any physical objects) and looking elsewhere or being distracted can be predictive of one's MA. It suggests that educators can use such behavioural observations as a quick gauge of whether a student is at risk of related socio-emotional struggles such as experiencing anxiety when learning math. While current findings might be limited in their generalisability to traditional educational settings and western countries, future studies should look to investigating the variables in other contexts.